\documentclass{article}
\usepackage[utf8]{inputenc}

\newtheorem{lemma}{Lemma}[section]

\newtheorem{theorem}[lemma]{Theorem}
\newtheorem{proposition}[lemma]{Proposition}
\newtheorem{corollary}[lemma]{Corollary}
\newtheorem{example}[lemma]{Example}

\newtheorem{note}[lemma]{Remark}

\title{Almost $n$-ary and almost $n$-aritizable theories\footnote{The work of the author was carried out
in the framework of the State Contract of the Sobolev Institute of
Mathematics, Project No.~FWNF-2022-0012 (Section 2), and of
Russian Scientific Foundation, Project No.~22-21-00044 (Sections
3, 4).}}
\author{S.V. Sudoplatov}
\date{}

\begin{document}

\maketitle
\begin{abstract}
We study possibilities for almost $n$-ary and $n$-aritizable
theories. Their dynamics both in general case, for
$\omega$-categorical theories, and with respect to operations for
theories are described.
\end{abstract}


\bigskip
We continue to study arities of theories and of their expansions
\cite{aafot}. In the present paper we introduce natural notions of
almost $n$-ary and almost $n$-aritizable theories, and describe
their dynamics both in general case, for $\omega$-categorical
theories, and with respect to operations for theories.

\section{Preliminaries}

Recall a series of notions related to arities and aritizabilities
of theories.

{\bf Definition} \cite{ZPS}. A theory $T$ is said to be
\emph{$\Delta$-based}\index{Theory!$\Delta$-based}, where $\Delta$
is some set of~formulae without parameters, if any formula of $T$
is equivalent in $T$ to a~Boolean combination of formulae
in~$\Delta$.

For $\Delta$-based theories $T$, it is also said that $T$ has {\em
quantifier elimination}\index{Elimination of quantifiers} or {\em
quantifier reduction}\index{Reduction of quantifiers} up to
$\Delta$.

\medskip
{\bf Definition} \cite{ZPS, CCMCT}. {\rm Let $\Delta$ be a set of
formulae of a theory $T$, and $p(\bar{x})$ a type of $T$ lying in
$S(T)$. The type $p(\bar{x})$ is said to be
\emph{$\Delta$-based}\index{Type!$\Delta$-based} if $p(\bar{x})$
is isolated by a set of formulas $\varphi^\delta\in p$, where
$\varphi\in\Delta$, $\delta\in\{0,1\}$.}

\medskip
The following lemma, being a corollary of Compactness Theorem,
noticed in \cite{ZPS}.

\medskip
\begin{lemma}\label{lem01} A  theory  $T$  is
$\Delta$-based  if and only if,  for  any  tuple~$\bar{a}$ of any
{\rm (}some{\rm )} weakly saturated model of $T$, the type ${\rm
tp}(\bar{a})$ is $\Delta$-based.
\end{lemma}

\medskip
{\bf Definition} \cite{aafot}. An elementary theory $T$ is called
{\em unary}, or {\em $1$-ary}, if any $T$-formula
$\varphi(\overline{x})$ is $T$-equivalent to a Boolean combination
of $T$-formulas, each of which is of one free variable, and of
formulas of form $x\approx y$.

For a natural number $n\geq 1$, a formula $\varphi(\overline{x})$
of a theory $T$ is called {\em $n$-ary}, or an {\em $n$-formula},
if $\varphi(\overline{x})$ is $T$-equivalent to a Boolean
combination of $T$-formulas, each of which is of $n$ free
variables.

For a natural number $n\geq 2$, an elementary theory $T$ is called
{\em $n$-ary}, or an {\em $n$-theory}, if any $T$-formula
$\varphi(\overline{x})$ is $n$-ary.

A theory $T$ is called {\em binary} if $T$ is $2$-ary, it is
called {\em ternary} if $T$ is $3$-ary, etc.

We will admit the case $n=0$ for $n$-formulae
$\varphi(\overline{x})$. In such a case $\varphi(\overline{x})$ is
just $T$-equivalent to a sentence
$\forall\overline{x}\varphi(\overline{x})$.

If $T$ is a theory such that $T$ is $n$-ary and not $(n-1)$-ary
then the value $n$ is called the arity of $T$ and it is denoted by
${\rm ar}(T)$. If $T$ does not have any arity we put ${\rm
ar}(T)=\infty$.

Similarly, for a formula $\varphi$ of a theory $T$ we denote by
${\rm ar}_T(\varphi)$ the natural value $n$ if $\varphi$ is
$n$-ary and not $(n-1)$-ary. If $\varphi$ does not any arity we
put ${\rm ar}_T(\varphi)=\infty$. If a theory $T$ is fixed we
write ${\rm ar}(\varphi)$ instead of ${\rm ar}_T(\varphi)$.

\medskip
The following example illustrates the notions above, and it will
be used below.

\begin{example}\label{excircular} {\rm Recall \cite{KulMac, Kulp3, Kulp4}
that a {\em circular}, or {\em cyclic}  order relation is
described by a ternary relation $K_3$ satisfying the following conditions:\\
 (co1) $\forall x\forall y \forall z (K_3(x,y,z)\to K_3(y,z,x));$\\
 (co2) $\forall x\forall y \forall z (K_3(x,y,z)\land K_3(y,x,z)
 \leftrightarrow x=y \lor y=z \lor z=x);$\\
 (co3) $\forall x\forall y \forall z(K_3(x,y,z)\to \forall t[K_3(x,y,t)
 \lor K_3(t,y,z)]);$\\
 (co4) $\forall x\forall y \forall z (K_3(x,y,z)\lor K_3(y,x,z)).$

Clearly, ${\rm ar}(K_3(x,y,z))=3$ if the relation has at least
three element domain. Hence, theories with infinite circular order
relations are at least $3$-ary.

The following generalization of circular order produces a {\em
$n$-ball}, or {\em $n$-spherical}, or {\em $n$-circular} order
relation, for $n\geq 4$, which is described by a $n$-ary relation
$K_n$ satisfying the following conditions:\\
 (nbo1) $\forall x_1,\ldots,x_n (K_n(x_1,x_2,\ldots,x_n)\to K_n(x_2,\ldots,x_n,x_1));$\\
 (nbo2) $\forall x_1,\ldots,x_n \biggm(K_n(x_1,\ldots,x_i,x_{i+1},\ldots,x_n)\land$  $$\land K_n(x_1,\ldots,x_{i+1},x_{i},\ldots,x_n)
 \leftrightarrow\bigvee\limits_{i=1}^{n-1} x_i=x_{i+1}\biggm);$$\\
 (nbo3) $\forall x_1,\ldots,x_n(K_n(x_1,\ldots,x_n)\to \forall t[K_n(x_1,\ldots,x_{n-1},t)
 \lor K_n(t,x_2,\ldots,x_n)]);$\\
 (nbo4) $\forall x_1,\ldots,x_n (K_n(x_1,\ldots,x_i,x_{i+1},\ldots,x_n)\lor K_n(x_1,\ldots,x_{i+1},x_{i},\ldots,x_n)),$ $i<n.$

Clearly, ${\rm ar}(K_n(x_1,\ldots,x_n))=n$ if the relation has at
least $n$-element domain. Thus, theories with infinite $n$-ball
order relations are at least $n$-ary.}
\end{example}

{\bf Definition} \cite{aafot}. A $T$-formula
$\varphi(\overline{x})$ is called {\em $n$-expansible}, or {\em
$n$-arizable}, or {\em $n$-aritizable}, if $T$ has an expansion
$T'$ such that $\varphi(\overline{x})$ is $T'$-equivalent to a
Boolean combination of $T'$-formulas with $n$ free variables.

A theory $T$ is called {\em $n$-expansible}, or {\em
$n$-arizable}, or {\em $n$-aritizable}, if there is an $n$-ary
expansion $T'$ of $T$.

A theory $T$ is called {\em arizable} or {\em aritizable}, if $T$
is $n$-aritizable for some $n$.

A $1$-aritizable theory is called {\em unary-able}, or {\em
unary-tizable}. A $2$-aritizable theory is called {\em
binary-tizable} or {\em binarizable}, a $3$-aritizable theory is
called {\em ternary-tizable} or {\em ternarizable}, etc.

\medskip
{\bf Definition.} \cite{Wo} The {\em disjoint
union}\index{Disjoint union!of structures}
$\bigsqcup\limits_{n\in\omega}\mathcal{
M}_n$\index{$\bigsqcup\limits_{n\in\omega}\mathcal{ M}_n$} of
pairwise disjoint structures $\mathcal{ M}_n$ for pairwise
disjoint predicate languages $\Sigma_n$, $n\in\omega$, is the
structure of language
$\bigcup\limits_{n\in\omega}\Sigma_n\cup\{P^{(1)}_n\mid
n\in\omega\}$ with the universe $\bigsqcup\limits_{n\in\omega}
M_n$, $P_n=M_n$, and interpretations of predicate symbols in
$\Sigma_n$ coinciding with their interpretations in $\mathcal{
M}_n$, $n\in\omega$. The {\em disjoint union of
theories}\index{Disjoint union!of theories} $T_n$ for pairwise
disjoint languages $\Sigma_n$ accordingly, $n\in\omega$, is the
theory
$$\bigsqcup\limits_{n\in\omega}T_n\rightleftharpoons{\rm Th}\left(\bigsqcup\limits_{n\in\omega}\mathcal{ M}_n\right),$$
where\index{$\bigsqcup\limits_{n\in\omega}T_n$} $\mathcal{
M}_n\models T_n$, $n\in\omega$. Taking empty sets instead of some
structures $\mathcal{ M}_k$ we obtain disjoint unions of finitely
many structures and theories. In particular, we have the disjoint
unions $\mathcal{ M}_0\sqcup\ldots\sqcup\mathcal{ M}_n$ and their
theories $T_0\sqcup\ldots\sqcup T_n$.

\begin{theorem}\label{thdu} {\rm \cite{aafot}.}
$1.$ For any theories $T_m$, $m\in\omega$, and their disjoint
union $\bigsqcup\limits_{m\in\omega}T_m$, all $T_m$ are
$n$-theories iff $\bigsqcup\limits_{m\in\omega}T_m$ is an
$n$-theory, moreover, $${\rm
ar}\left(\bigsqcup\limits_{m\in\omega}T_m\right)={\rm max}\{{\rm
ar}(T_m)\mid m\in\omega\}.$$

$2.$ For any theories $T_m$, $m\in\omega$, and their disjoint
union $\bigsqcup\limits_{m\in\omega}T_m$, all $T_m$ are
$n$-aritizable iff $\bigsqcup\limits_{m\in\omega}T_m$ is
$n$-aritizable.
\end{theorem}

{\bf Definition} \cite{EKS20}. Let $\mathcal{M}$ and $\mathcal{N}$
be structures of relational languages $\Sigma_\mathcal{M}$ and
$\Sigma_\mathcal{N}$ respectively. We define the {\em composition}
$\mathcal{M}[\mathcal{N}]$ of $\mathcal{M}$ and $\mathcal{N}$
satisfying the following conditions:

1)
$\Sigma_{\mathcal{M}[\mathcal{N}]}=\Sigma_\mathcal{M}\cup\Sigma_\mathcal{N}$;

2) $M[N]=M\times N$, where $M[N]$, $M$, $N$ are universes of
$\mathcal{M}[\mathcal{N}]$, $\mathcal{M}$, and $\mathcal{N}$
respectively;

3) if $R\in\Sigma_\mathcal{M}\setminus\Sigma_\mathcal{N}$,
$\mu(R)=n$, then $((a_1,b_1),\ldots,(a_n,b_n))\in
R_{\mathcal{M}[\mathcal{N}]}$ if and only if $(a_1,\ldots,a_n)\in
R_{\mathcal{M}}$;

4) if $R\in\Sigma_\mathcal{N}\setminus\Sigma_\mathcal{M}$,
$\mu(R)=n$, then $((a_1,b_1),\ldots,(a_n,b_n))\in
R_{\mathcal{M}[\mathcal{N}]}$ if and only if $a_1=\ldots =a_n$ and
$(b_1,\ldots,b_n)\in R_{\mathcal{N}}$;

5) if $R\in\Sigma_\mathcal{M}\cap\Sigma_\mathcal{N}$, $\mu(R)=n$,
then $((a_1,b_1),\ldots,(a_n,b_n))\in
R_{\mathcal{M}[\mathcal{N}]}$ if and only if $(a_1,\ldots,a_n)\in
R_{\mathcal{M}}$, or $a_1=\ldots =a_n$ and $(b_1,\ldots,b_n)\in
R_{\mathcal{N}}$.

The composition $\mathcal{M}[\mathcal{N}]$ is called {\em
$e$-definable}, or {\em {\rm equ}-definable}, if
$\mathcal{M}[\mathcal{N}]$ has an $\emptyset$-definable
equivalence relation $E$ whose $E$-classes are universes of the
copies of $\mathcal{N}$ forming $\mathcal{M}[\mathcal{N}]$. If the
equivalence relation $E$ is fixed, the $e$-definable composition
is called {\em $E$-definable}.

\medskip
Using a nice basedness of $E$-definable compositions $T_1[T_2]$
(see \cite{EKS20}) till the formulas of form $E(x,y)$ and
generating formulas for $T_1$ and $T_2$ we have the following:

\begin{theorem}\label{thcomp} {\rm \cite{aafot}.} $1.$  For any theories $T_1$ and $T_2$ and
their $E$-definable composition $T_1[T_2]$, $T_1$ and $T_2$ are
$n$-theories, for $n\geq 2$, iff $T_1[T_2]$ is an $n$-theory,
moreover, ${\rm ar}(T_1[T_2])={\rm max}\{{\rm ar}(T_1),{\rm
ar}(T_2)\}$, if models of $T_1$ and of $T_2$ have at least two
elements, and ${\rm ar}(T_1[T_2])={\rm max}\{{\rm ar}(T_1),{\rm
ar}(T_2),2\}$, if a model of $T_1$ or $T_2$ is a singleton.

$2.$ For any theories $T_1$ and $T_2$ and their $E$-definable
composition $T_1[T_2]$, $T_1$ and $T_2$ are $n$-aritizable iff
$T_1[T_2]$ is $n$-aritizable.
\end{theorem}

\section{Almost $n$-ary and $n$-aritizable theories, their dynamics}

{\bf Definition.} (Cf. \cite{Kulp3, Kulp4}) A theory $T$ is called
{\em almost $n$-ary} if there are finitely many formulae
$\varphi_1(\overline{x}),\ldots,\varphi_m(\overline{x})$ such that
each $T$-formula is $T$-equivalent to a Boolean combination of
$n$-formulae and formulae obtained by substitutions of free
variables in
$\varphi_1(\overline{x}),\ldots,\varphi_m(\overline{x})$.

In such a case we say that the formulae
$\varphi_1(\overline{x}),\ldots,\varphi_m(\overline{x})$ witness
that $T$ is almost $n$-ary.

Almost $1$-ary theories are called {\em almost unary}, almost
$2$-ary theories are called {\em almost binary}, almost $3$-ary
theories are called {\em almost ternary}, etc.

A theory $T$ is called {\em almost $n$-aritizable} if some
expansion $T'$ of $T$ is almost $n$-ary.

Almost $1$-aritizable theories are called {\em almost
unary-tizable}, almost $2$-aritizable theories are called {\em
almost binarizable}, almost $3$-aritizable theories are called
{\em almost ternarizable}, etc.

\medskip
The following properties are obvious.

\medskip 1. Any $n$-ary (respectively, $n$-aritizable) theory is
almost $n$-ary (almost $n$-aritizable).

\medskip
2. Any almost $n$-ary (respectively, $n$-aritizable) theory is
almost $k$-ary (almost $k$-aritizable) for any $k\geq n$.

\medskip
3. Any theory of a finite structure is almost unary.

\medskip
Families of weakly circularly minimal structures produce examples
of almost binary theories which are not binary \cite{Kulp3,
Kulp4}. Similarly natural generalizations of weakly circularly
minimal structures till $n$-circular orders give examples of
almost $(n-1)$-ary theories $T_n$ with ${\rm ar}(T_n)=n$, $n\geq
4$.

\medskip
Assuming that the witnessing set
$\{\varphi_1(\overline{x}),\ldots,\varphi_m(\overline{x})\}$ is
minimal for the almost $n$-ary theory $T$ we have either $m=0$ of
$l(\overline{x})>n$.

Thus we have two minimal characteristics witnessing the almost
$n$-arity of $T$: $m$ and $l(\overline{x})$. The pair
$(m,l(\overline{x}))$ is called the {\em degree} of the almost
$n$-arity of $T$, or the {\em {\rm aar}-degree} of $T$, denoted by
${\rm deg}_{\rm aar}(T)$. Here we assume that $n$ is minimal with
almost $n$-arity of $T$, this $n$ is denoted by ${\rm aar}(T)$.
Clearly, ${\rm aar}(T)\leq{\rm ar}(T)$, and if $m=0$, i.e.,
$n={\rm ar}(T)={\rm aar}(T)$ then it is supposed that
$l(\overline{x})=0$, too.

We have ${\rm aar}(T)\in\omega$ if and only if ${\rm
ar}(T)\in\omega$. So if ${\rm ar}(T)=\infty$ then it is natural to
put ${\rm aar}(T)=\infty$.

Besides, $n={\rm ar}(T)={\rm aar}(T)\in\omega$ means that the set
$\Delta_n(T)$ of $T$-formulae with $n$ free variables allows to
express all $\emptyset$-definable sets for $T$ by Boolean
combinations and taking any set $\Delta_k(T)$ of $T$-formulae with
$k<n$ free variables all $\emptyset$-definable sets for $T$ can be
expressed by Boolean combinations of formulae in $\Delta_k$ and
substitutions of infinitely many formulae $\varphi(\overline{x})$
only, where $l(\overline{x})=n$.

By the definition if ${\rm aar}(T)=n$ then
\begin{equation}\label{eq_aar}{\rm
deg}_{\rm ar}(T)\in\{(0,0)\}\cup\{(m,r)\mid
m\in\omega\setminus\{0\}, r\in\omega, r>n\}.
\end{equation}
The described pairs in the relation (\ref{eq_aar}) are called {\em
admissible}.

\begin{theorem}\label{th_aar_mn}
For any $m,n\in\omega\setminus\{0\}$ with $m\leq n$ there is a
theory $T_{mn}$ with ${\rm aar}(T{mn})=m$ and ${\rm
ar}(T_{mn})=n$.
\end{theorem}

Proof. We use disjoint unions of {\em dense} $n$-spherically
ordered theories $T_n$ with infinite orders, i.e., theories,
generated by $n$-spherical orders $K_n(x_1,x_2,x_3,\ldots,x_n)$
satisfying the axioms $$\forall x_1,x_2,\ldots
x_n\Bigg(K_n(x_1,x_2,x_3,\ldots,x_n)\wedge\bigwedge_{i\ne j}\neg
x_i\approx x_j\to$$
$$\to\exists y\Bigg(\bigwedge_{i\leq n}\neg
x_i\approx y\wedge K_n(x_1,y,x_3,\ldots,x_n)\Bigg)\Bigg).$$ For
$n=2$ we take the theory $T_2$ of dense linear order
$K_2(x_1,x_2)$ without endpoints, having ${\rm ar}(T_2)=2$, and
for $n=1$~--- the theory $T_1$ of the empty languages, having
${\rm ar}(T_1)=1$.

Similarly to dense linear orders and dense circular orders, dense
$n$-spherical orders produce quantifier eliminations with ${\rm
ar}(T_n)=n$, $n\in\omega\setminus\{0\}$.

Using Theorem \ref{thdu} we obtain $m={\rm ar}(T)={\rm aar}(T)$
taking $T_{mm}=\bigsqcup\limits_{r\in\omega}T^r_m$, where $T^r_m$
are copies of $T_m$ in disjoint languages $\{K^r_m\}$,
$r\in\omega$. Indeed, disjoint predicates $K^r_m$ producing ${\rm
ar}(T_m)=m$ witness that ${\rm ar}(T_{mm})=m$. Finitely many these
predicates can not define all definable sets for $T_{mm}$ since
there are infinitely many of them. Thus, ${\rm aar}(T_{mm})=m$,
too.

Now for any $m<n$ we form $T_{mn}=T_n\sqcup
\bigsqcup\limits_{r\in\omega}T^r_m$. By $T_n$ in $T_{mn}$ and
$m<n$ we have ${\rm ar}(T_{mn})=n$. And ${\rm aar}(T_{mn})=m$
since there are infinitely many disjoint predicates of arity $m$,
as required.

\medskip
The following theorem shows that all admissible pairs are
realized.

\begin{theorem}\label{th_admis}
For any admissible pair $(m,r)$ and $n\in\omega\setminus\{0,1\}$
there is a theory $T$ with ${\rm aar}(T)=n$ and ${\rm deg}_{\rm
ar}(T)=(m,r)$.
\end{theorem}

Proof. The admissible pair $(0,0)$ with ${\rm ar}(T)=n$ is
realized by Theorem \ref{th_aar_mn}. Now for an admissible pair
$(m,r)\ne(0,0)$ we can take a disjoint union $T$ of countably many
dense $n$-spherically ordered theories and of $m$ dense
$r$-spherically ordered theories. Using arguments for Theorem
\ref{th_aar_mn} we obtain ${\rm aar}(T)=n$ and ${\rm deg}_{\rm
ar}(T)=(m,r)$, as required.

\begin{proposition}\label{praar1}
Any almost $n$-ary theory $T$ is $k$-ary for some $k\geq n$.
\end{proposition}

Proof. Let $T$ be an almost $n$-ary theory witnessed by the
formulae $\varphi_1(\overline{x}),\ldots,$
$\varphi_m(\overline{x})$. Then taking the set $\Delta_k$ of all
formulae with $k={\rm max}\{n,l(\overline{x})\}$ we observe, using
$\varphi_1(\overline{x}),\ldots,\varphi_m(\overline{x})$, that $T$
is $\Delta_k$-based, i.e., $T$ is $k$-ary, as required.

\begin{corollary}\label{coraar1}
Any theory $T$ is $n$-ary for some $n$ iff $T$ is almost $m$-ary
for some $m$.
\end{corollary}

\begin{corollary}\label{coraar2}
Any theory $T$ is $n$-aritizable for some $n$ iff $T$ is almost
$m$-aritizable for some $m$.
\end{corollary}

\section{$\omega$-categorical almost $n$-ary and $n$-aritizable theories}

\begin{proposition}\label{praar_cat1}
If $T$ is an almost $n$-ary $\omega$-categorical theory, for some
$n$, then $T$ is almost $k$-ary for any
$k\in\omega\setminus\{0\}$, i.e., ${\rm aar}(T)=1$.
\end{proposition}

Proof. If $k\geq n$ then $T$ is almost $m$-ary, as noticed above.
If $k<n$ then we collect in a set $Z$ all formulae
$\varphi_1(\overline{x}),\ldots,\varphi_m(\overline{x})$
witnessing the almost $n$-arity of $T$ and, by Ryll-Nardzewski
Theorem, all non-equivalent formulae
$\psi_1(\overline{y}),\ldots,\psi_r(\overline{y})$ with
$l(\overline{y})=n$. Clearly, the set $Z$ witnesses that $T$ is
almost $k$-ary. Taking $k=1$ we obtain ${\rm aar}(T)=1$, as
required.

\begin{corollary}\label{coraar3}
For any $\omega$-categorical theory $T$ either ${\rm aar}(T)=1$
with ${\rm ar}(T)\in\omega$, or ${\rm aar}(T)=\infty$ with ${\rm
ar}(T)=\infty$.
\end{corollary}

The following example illustrates Corollary \ref{coraar3}.

\begin{example}
{\rm Taking a dense linear order $K_2$ without endpoint we can
step-by-step extend it to a chain of dense $n$-spherical orders
$K_n$, $n\geq 2$, in the following way.

We put $(a,b,c)\in K_3$ if $(a,b)\in K_2$ and $(b,c)\in K_2$, or
$(b,c)\in K_2$ and $(c,a)\in K_2$, or $(c,a)\in K_2$ and $(a,b)\in
K_2$. If $K_n$, $n\geq 3$, is defined then we put
$(a_1,\ldots,a_{n+1})\in K_{n+1}$ if $(a_1,\ldots,a_n)\in K_n$ and
$(a_2,\ldots,a_{n+1})\in K_n$, or $a_2,\ldots,a_{n+1}\in K_n$ and
$(a_3,\ldots,a_{n+1},a_1)\in K_n$, or $(a_3,\ldots,a_{n+1},a_1)\in
K_n$ and $(a_4,\ldots,a_{n+1},a_1,a_2)\in K_n$. The obtained
structure $\mathcal{M}_\infty$ in the language
$\Sigma_\infty=\{K_n\mid n\in\omega\setminus\{0,1\}\}$ has an
$\omega$-categorical theory $T_\infty$ with ${\rm
aar}(T_\infty)={\rm ar}(T_\infty)=\infty$, since each new $K_n$
increases the arity. The same characteristics have restrictions of
$T_\infty$ to any infinite sublanguages.

At the same time each restriction $\mathcal{M}$ of
$\mathcal{M}_\infty$ to a finite nonempty sublanguage
$\{K_{n_1},\ldots,$ $K_{n_m}\}$ produces a theory $T$ with ${\rm
aar}(T)=1$ with ${\rm ar}(T)={\rm max}\{n_1,\ldots,n_m\}$.}
\end{example}

By Corollary \ref{coraar3} and the definition of almost
aritizability we immediately have:

\begin{corollary}\label{coraar4}
Any restriction $T$ of $n$-ary $\omega$-categorical theory $T'$ is
almost unary-tizable.
\end{corollary}

\section{Operations for almost $n$-ary and $n$-aritizable theories}

In this section we consider links for arities of theories with
respect to disjoint unions of theories and $E$-definable
compositions of theories.

\begin{theorem}\label{thduaa}
$1.$ For any theories $T_1$, $T_2$ and their disjoint union
$T_1\sqcup T_2$, both $T_1$ and $T_2$ are almost $n$-ary iff
$T_1\sqcup T_2$ is almost $n$-ary, moreover, ${\rm aar}(T_1\sqcup
T_2)={\rm max}\{{\rm aar}(T_1),{\rm aar}(T_2)\}$.

$2.$ For any theories $T_1$, $T_2$ and their disjoint union
$T_1\sqcup T_2$, both $T_1$ and $T_2$ are almost $n$-aritizable
iff $T_1\sqcup T_2$ is almost is $n$-aritizable.
\end{theorem}

Proof. 1. Let $\Phi_1$ and $\Phi_2$ be finite sets of formulas
witnessing that $T_1$ and $T_2$ are almost $n$-ary, respectively.
Using the definition of disjoint union and Theorem \ref{thdu} we
obtain that the finite set $\Phi_1\cup \Phi_2$ witnesses that
$T_1\sqcup T_2$ is almost $n$-ary. Conversely, if a finite set
$\Phi$ of formulas witnesses that $T_1\sqcup T_2$ is almost
$n$-ary then $Phi$ witnesses that both $T_1$ and $T_2$ are almost
$n$-ary.

The equality ${\rm aar}(T_1\sqcup T_2)={\rm max}\{{\rm
aar}(T_1),{\rm aar}(T_2)\}$ follows from the definition of
disjoint union since if ${\rm aar}(T_1\sqcup T_2)=k$ then the
maximal value of ${\rm aar}(T_1)$ and ${\rm aar}(T_2)$ is
responsible for this equality.

Item 2 follows from Item 1 since expansions of $T_1$ and $T_2$
correspond expansions of $T_1\sqcup T_2$: some expansions of
$T'_1$ and $T'_2$ of $T_1$ and $T_2$, respectively, are almost
$n$-ary iff $T'_1\sqcup T'_2$ produces an almost $n$-ary expansion
of $T_1\sqcup T_2$, as required.

\medskip
Using induction we obtain the following:

\begin{corollary}\label{coduaa}
$1.$ For any theories $T_1, T_2,\ldots,T_m$ and their disjoint
union $\bigsqcup\limits_{i=1}^m T_i$, all $T_1, T_2,\ldots,T_m$
are almost $n$-ary iff $\bigsqcup\limits_{i=1}^m T_i$ is almost
$n$-ary, moreover, $${\rm aar}\left(\bigsqcup\limits_{i=1}^m
T_i\right)={\rm max}\{{\rm aar}(T_i)\mid i\leq m\}.$$

$2.$ For any theories $T_1, T_2,$ $\ldots,T_m$ and their disjoint
union $\bigsqcup\limits_{i=1}^m T_i$, all $T_1, T_2,$ $\ldots,T_m$
are almost $n$-aritizable iff $\bigsqcup\limits_{i=1}^m T_i$ is
almost $n$-aritizable.
\end{corollary}

\begin{note}{\rm Both almost $n$-arity and almost
$n$-aritizability can fail taking disjoint unions of infinitely
many theories $T_i$, $i\in I$. Indeed, each theory $T_i$ can have
its own finite set $\Phi_i$ of formulas witnessing the almost
$n$-arity/$n$-aritizability, say in disjoint languages, whereas
finite unions $\bigcup\Phi_i$ can not witness the almost
$n$-arity/$n$-aritizability for $\bigsqcup\limits_{i\in I}T_i$.}
\end{note}

Generalizing Theorem \ref{thcomp} we obtain:

\begin{theorem}\label{thcompaa} $1.$ For any theories $T_1$ and $T_2$ and
their $E$-definable composition $T_1[T_2]$, $T_1$ and $T_2$ are
almost $n$-ary iff $T_1[T_2]$ is almost $n$-ary, moreover, ${\rm
aar}(T_1[T_2])={\rm max}\{{\rm aar}(T_1),{\rm aar}(T_2),2\}$, if
models of $T_1$ and of $T_2$ have at least two elements, and ${\rm
aar}(T_1[T_2])={\rm max}\{{\rm aar}(T_1),{\rm aar}(T_2)\}$, if a
model of $T_1$ or $T_2$ is a singleton.

$2.$ For any theories $T_1$ and $T_2$ and their $E$-definable
composition $T_1[T_2]$, $T_1$ and $T_2$ are almost $n$-aritizable
iff $T_1[T_2]$ is almost $n$-aritizable.
\end{theorem}

Proof. 1. Let $T_i$ be $\Delta_i$-based for $i=1,2$. Since
$T_1[T_2]$ is $E$-definable it is $\Delta$-based, where $\Delta$
consists of formulae in $\Delta_1\cup\Delta_2$ and $E(x,y)$
\cite{EKS20}. Now assuming that $T_1$ and $T_2$ are almost $n$-ary
we can choose $\Delta_i$ consisting of $n$-formulae and finitely
many formulae forming $\Phi_i$, $i=1,2$. Hence $T_1[T_2]$ is
almost $n$-ary.

Conversely, if $T_1[T_2]$ is almost $n$-ary and it is witnessed by
a set $\Phi$ of formulae then $\Phi$ witnesses that $T_1$ and
$T_2$ are almost $n$-ary.

If models of $T_1$ and of $T_2$ have at least two elements then
$T_1[T_2]$ is at least binary that witnessed by the formula
$E(x,y)$. Thus since $T_1[T_2]$ is $\Delta$-based we have ${\rm
aar}(T_1[T_2])={\rm max}\{{\rm aar}(T_1),{\rm aar}(T_2),2\}$. If
$T_1$ or $T_2$ is a theory of singleton then $T_1[T_2]$ is
$(\Delta_1\cup\Delta_2)$-based implying ${\rm aar}(T_1[T_2])={\rm
max}\{{\rm aar}(T_1),{\rm aar}(T_2)\}$.

Item 2 follows from Item 1 repeating the arguments for Item 2 of
Theorem \ref{thduaa}.

\medskip
Theorem \ref{thcompaa} immediately implies

\begin{corollary}\label{cocompaa}
Any composition of finitely many almost $n$-ary {\rm (}almost
$n$-ari\-ti\-z\-a\-b\-le{\rm )} theories, for $n\geq 2$, is again
many almost $n$-ary {\rm (}almost $n$-aritizable{\rm )}.
\end{corollary}

\section{Conclusion}

We considered possibilities for almost arities and almost
aritizabilities of theories and their dynamics both in general
case, for $\omega$-categorical theories, and with respect to
operations for theories. It would be interesting to describe
values of almost arities and almost aritizabilities for natural
classes of theories.

\noindent Sobolev Institute of Mathematics, \\ 4, Acad. Koptyug
avenue, Novosibirsk, 630090, Russia; \\ Novosibirsk State
Technical
University, \\ 20, K.Marx avenue, Novosibirsk, 630073, Russia; \\
Novosibirsk State University, \\ 1, Pirogova street, Novosibirsk,
630090, Russia

\medskip\noindent
e-mail: sudoplat@math.nsc.ru

\end{document}